\documentclass[12pt,reqno]{amsart}%
\usepackage{amsfonts,anysize}
\usepackage{mathpazo}
\usepackage[colorlinks,urlcolor=blue,citecolor=blue,linkcolor=blue]{hyperref}
\usepackage{amsmath}
\usepackage{amsfonts}
\usepackage{amssymb}
\usepackage{graphicx}%
\providecommand{\U}[1]{\protect\rule{.1in}{.1in}}
\marginsize{2cm}{2cm}{1cm}{1cm}

\newtheorem{thm}{Theorem}
\newtheorem{cor}[thm]{Corollary}

\newtheorem{lem}[thm]{Lemma}

\theoremstyle{definition}

\theoremstyle{remark}
\newtheorem{rem}[thm]{Remark}
\newtheorem*{pf}{Proof}
\allowdisplaybreaks
\begin{document}

\title[surjectivity of smooth maps into Euclidean spaces]{On surjectivity of smooth maps into Euclidean spaces and the fundamental
theorem of algebra}
\author{Peng Liu
\and Shibo Liu}
\dedicatory{Department of Mathematics, Xiamen University\\
Xiamen 361005, China}
\thanks{Emails: \texttt{liupeng1729@qq.com} (P. Liu), \texttt{liusb@xmu.edu.cn} (S.B. Liu)}

\begin{abstract}
In this note we obtain the surjectivity of smooth maps into Euclidean spaces
under mild conditions. As application we give a new proof of the Fundamental
Theorem of Algebra. We also observe that any $C^1$-map from a compact manifold into Euclidean space with dimension $n\ge2$ has infinitely many critical points.
\end{abstract}
\maketitle

Let $f:\mathbb{R}^{n}\rightarrow\mathbb{R}^{n}$ be a $C^{1}$-map such that for
all $x\in\mathbb{R}^{n}$ we have $\det Df(x)\neq0$. Suppose moreover%
\begin{equation}
\lim_{\left\vert x\right\vert \rightarrow\infty}\left\vert f(x)\right\vert
=+\infty\text{,} \label{e1}%
\end{equation}
then it is well known that $f(\mathbb{R}^{n})=\mathbb{R}^{n}$ (namely $f$ is
surjective), see e.g. \cite[Page 24]{MR787404}. In this note, we will prove
the following generalisation of this result. As we shall see, from this
result, we can easily obtain the fundamental theorem of algebra.

Let $U$ be an open subset of $\mathbb{R}^{m}$ and $f:U\rightarrow
\mathbb{R}^{n}$ be a $C^{1}$-map, $a\in U$. If $\operatorname*{rank}Df(a)<n$,
that is the differential map $Df(a):\mathbb{R}^{m}\rightarrow\mathbb{R}^{n}$
is not surjective, then we say that $a$ is a critical point of $f$.

\begin{thm}
\label{t1}Let $f:\mathbb{R}^{m}\rightarrow\mathbb{R}^{n}$ be a $C^{1}$-map with $n\ge2$. If
$f$ has only finitely many critical points and $f(\mathbb{R}^{m})$ is a closed
subset of $\mathbb{R}^{n}$, then $f(\mathbb{R}^{m})=\mathbb{R}^{n}$.
\end{thm}

\begin{rem}
If $m=n$, then $\det Df(x)\neq0$ is exactly that $x$ is not a critical point
of $f$. If the condition \eqref{e1} is satisfied, then $f(\mathbb{R}^{n})$ is
closed. Therefore the classical result mentioned at the beginning is a
corollary of our Theorem \ref{t1}.
\end{rem}

To prove this theorem, we need the following lemma.

\begin{lem}
\label{l1}Let  $n\ge2$, $A$ be a nonempty open subset of $\mathbb{R}^{n}$. If there
exists $k$ points $p_{1},\ldots,p_{k}$ such that $A\cup\left\{  p_{i}\right\}
_{i=1}^{k}$ is closed, then $\bar{A}=\mathbb{R}^{n}$.
\end{lem}

\begin{pf}
Because $A\cup\left\{  p_{i}\right\}  _{i=1}^{k}$ is closed, we have%
\[
A\cup\left\{  p_{i}\right\}  _{i=1}^{k}=\overline{A\cup\left\{  p_{i}\right\}
_{i=1}^{k}}=\bar{A}\cup\left\{  p_{i}\right\}  _{i=1}^{k}=A\cup\partial
A\cup\left\{  p_{i}\right\}  _{i=1}^{k}\text{.}%
\]
Since $A$ is open, $A\cap\partial A=\emptyset$. It follows that
\begin{equation}
\partial A\subset\left\{  p_{i}\right\}  _{i=1}^{k}\text{.} \label{e3}%
\end{equation}

If $\bar{A}\neq\mathbb{R}^{n}$, we can choose $b\in\mathbb{R}^{n}%
\backslash\bar{A}$ and $a\in A$. Since both $A$ and $\mathbb{R}^{n}%
\backslash\bar{A}$ are open, there exists $\varepsilon>0$ such that
\[
B_{\varepsilon}(a)\subset A\text{,\qquad}B_{\varepsilon}(b)\subset
\mathbb{R}^{n}\backslash\bar{A}\text{.}%
\]
This means that $b$ is not a boundary point of $A$.

Since  $n\ge2$, we can take a segment $\ell\subset B_{\varepsilon/2}(a)$ such that $\ell$ is not
parallel to the vector $b-a$. Of course we can choose $k+1$ different points
$\left\{  x_{i}\right\}  _{i=1}^{k+1}\subset\ell$.

Consider the $k+1$ segments $I_{i}$ with end points $b$ and $x_{i}$. Because
each such segment not only contains points (those near $x_{i}$) in $A$, but
also contains points (those near $b$) in $\mathbb{R}^{n}\backslash A$, it is
easy to see that for each $i=1,\ldots,k+1$, there exists at least one boundary
point of $A$ in $I_{i}$. Since for $i\neq j$ we have $I_{i}\cap I_{j}=\left\{
b\right\}  $ and $b\notin\partial A$, we conclude that $\partial A$ contains
at lease $k+1$ points. This contradicts with \eqref{e3}.
\end{pf}

\begin{pf}
[Proof of Theorem \ref{t1}]Let $K$ be the set of critical points, then $K$ is
a finite set. Because $\mathbb{R}^{m}\backslash K$ is open and
$\operatorname*{rank}Df(x)=n$ for $x\in\mathbb{R}^{m}\backslash K$, using the
Inverse Function Theorem it is well known that $A=f(\mathbb{R}^{m}\backslash
K)$ is an open subset of $\mathbb{R}^{n}$. By the assumption,%
\[
A\cup f(K)=f(\mathbb{R}^{m}\backslash K)\cup f(K)=f(\mathbb{R}^{m})
\]
is closed. Since $K$ is finite, its image $f(K)$ is also finite. Applying
Lemma \ref{l1}, we conclude that $\bar{A}=\mathbb{R}^{n}$.

Using the assumption that $f(\mathbb{R}^{m})$ is closed again, we deduce%
\[
f(\mathbb{R}^{m})=\overline{f(\mathbb{R}^{m})}\supset\overline{f(\mathbb{R}%
^{m}\backslash K)}=\bar{A}=\mathbb{R}^{n}\text{.}%
\]

\end{pf}

Checking the above proof of Theorem \ref{t1}, we find that the domain of our
map $f$ can be replaced by an $m$-dimensional $C^{1}$-manifold $M$. Therefore
we have the following corollary, whose proof is omitted.

\begin{cor}
Let $M$ be an $m$-dimensional $C^{1}$-manifold (without boundary), $f:M\rightarrow\mathbb{R}^{n}$
be a $C^{1}$-map  with $n\ge2$. If $f$ has only finitely many critical points and $f(M)$
is a closed subset of $\mathbb{R}^{n}$, then $f(M)=\mathbb{R}^{n}$.
\end{cor}

From this corollary, we can further state the following result.

\begin{cor}
Let $M$ be a compact $m$-dimensional $C^{1}$-manifold (without boundary), $f:M\rightarrow\mathbb{R}^{n}$
be a $C^{1}$-map with $n\ge2$. Then $f$ has infinitely many critical points.
\end{cor}

Now, we consider polynomials%
\[
p(z)=z^{n}+a_{1}z^{n-1}+\cdots+a_{n}%
\]
of degree $n\geq1$, where the coefficients $a_{i}\in\mathbb{C}$. As corollary
of our theorem we will give a new proof of the following classical theorem.

\begin{thm}
[Fundamental Theorem of Algebra]If $p(z)$ is a polynomial of degree $n\geq1$,
with complex coefficients $a_{i}$, then there exists $\xi\in\mathbb{C}$ such
that $p(\xi)=0$.
\end{thm}

\begin{pf}
Write $z=x+\mathrm{i}y$, and%
\[
p(z)=u(x,y)+\mathrm{i}v(x,y)\text{,}%
\]
we know that the complex derivative of $p$ is given by%
\[
p^{\prime}(z)=u_{x}+\mathrm{i}v_{x}\text{.}%
\]
If we consider our polynomial as a map $p:\mathbb{R}^{2}\rightarrow
\mathbb{R}^{2}$,%
\[
p(x,y)=\left(  u(x,y),v(x,y)\right)  \text{,}%
\]
then using the Cauchy-Riemann equations ($u_{x}=v_{y}$, $u_{y}=-v_{x}$), we
have%
\[
\det Dp(x,y)=\det\left(
\begin{array}
[c]{cc}%
u_{x} & u_{y}\\
v_{x} & v_{y}%
\end{array}
\right)  =\det\left(
\begin{array}
[c]{cr}%
u_{x} & -v_{x}\\
v_{x} & u_{x}%
\end{array}
\right)  =u_{x}^{2}+v_{x}^{2}\text{.}%
\]
Therefore, $p^{\prime}(z)=0$ if and only if $\det Dp(x,y)=0$. Because
$p^{\prime}(z)$ is a polynomial of degree $n-1$, it has at most $n-1$ zero
points. Thus, the map $p:\mathbb{R}^{2}\rightarrow\mathbb{R}^{2}$ has at most
$n-1$ critical points.

On the other hand, it is obvious that%
\[
\lim_{\left\vert \left(  x,y\right)  \right\vert \rightarrow\infty}\left\vert
p(x,y)\right\vert =+\infty\text{,}%
\]
which implies that $p(\mathbb{R}^{2})$ is closed. Applying Theorem \ref{t1} we
deduce $p(\mathbb{R}^{2})=\mathbb{R}^{2}$. In particular, $\left(  0,0\right)
\in p(\mathbb{R}^{2})$, this is equivalent to the existence of $\xi
\in\mathbb{C}$ such that $p(\xi)=0$.
\end{pf}

\begin{rem}
Although the Fundamental Theorem of Algebra has been proved for more than 200
years, new proofs of this theorem keep emerging even in the last decades, see e.g. \cite{MR2663254,MR1543753}. The spirit of our proof is somewhat similar to \cite{MR1543753}. In \cite{MR1543753}, the concepts of open (and closed) subsets of a topological subspace of $\mathbb{C}$ and the connectedness of such a subspace, are employed. In this sense, our proof is more elementary.
\end{rem}

\end{document}